\documentclass[12pt,a4paper]{article}
\usepackage{amssymb,amsmath}
\usepackage{hyperref}
\usepackage{amsthm}
\usepackage{ulem}
\newtheorem*{theorem}{Theorem}
\usepackage{chngcntr}
\usepackage{color}
\usepackage{graphicx}
\usepackage[dvipsnames]{xcolor}
\usepackage[symbol]{footmisc}
\usepackage[english,ngerman]{babel}
\usepackage[top=2.0cm,bottom=2.0cm,left=2.0cm,right=2.0cm]{geometry}

\begin{document}
\begin{center}
\LARGE{\textbf{
A rational series expansion of the Arctangent and a method for generating BBP series for Machin-Type formulas
}}
\[\]
Eine rationale Reihendarstellung für den Arkustangens und ein Verfahren zur Konstruktion von BBP Reihen für Formeln vom Machin'schen Typ
\[\]
\[\]
\large
{Alois Schiessl}
\footnote[1]{\,\;E-Mail: \texttt{aloisschiessl@posteo.de}}
\centerline{ }
\centerline{\today}
\end{center}
\selectlanguage{english}
\[\]
\[\]
\[\]
\begin{abstract}
In this work we derive a generalized series expansion of the real Arctangent, which is parametrized by by a natural number $M$. Each $M\in \mathbb{N}$ yields a new series. The expansion enables the construction of BPP formulas for the Arctangent. BBP formulas are of interest because, for fixed base $b$, the $k^{th}$ $b$-ary digit of a number with a BBP formula can be found without computing prior digits. The case $x=1$ leads to binary digit extraction formulas for $\frac{\pi}{4}$.
\[\]
\[\]
\[\textbf{Zusammenfassung}\]
In dieser Veröffentlichung geben wir für die reelle $\arctan\left(x\right)$ Funktion eine Reihendarstellung Darstellung an, die durch eine natürliche Zahl $M$ parametrisiert ist. Jedes $M\in \mathbb{N}$ führt zu einer neuen Reihe. Damit können BBP-Reihen für den Arkustangens konstruiert werden. BBP-Reihen sind deshalb von Interesse, weil mit ihnen die sich an der Stelle $k$ befindliche Zahl berechnet werden kann, ohne dass die Kenntnis der vorherigen Stellen nötig wäre. Im Spezialfall $x=1$ ergeben sich binäre BPP-Formeln für $\frac{\pi}{4}$.
\[\text{Deutsche Version ab Seite 12}\]
\end{abstract}
\selectlanguage{english}
\section{Introduction and main result}
In this paper, we present a novel representation of the real $\arctan\left(x\right)$ function, parameterized by $M\in \mathbb{N}$. Each natural number $M$ yields a new series. This allows for the construction of efficient BBP series for Arctangents. BBP formulas are of interest because, for fixed $b$, the $k^{th}$ $b$-ary digit of a number with a BBP formula can be found without computing prior digits.
$\\ $
$\\ $
Now we're going to formulate our statement
\begin{theorem}
$ \\ $
Let $1-\sqrt{2}< x < 1+\sqrt{2}$ and $M\in\mathbb N$. Then, the following statement holds.
\begin{align*}
&\arctan\left(x\right)=\sum_{k=0}^{\infty}
\left(\left(\frac{-4\,x^4}{\left(x+1\right)^4}\right)^M\right)^{k}
\left(\sum_{\ell=0}^{M-1}\left(\frac{-4\,x^4}{\left(x+1\right)^4}\right)^{\ell}
\right.\times
\\
&\Biggl.
\left(
 \frac{1}{{4\,k\,M}+4\,\ell+1}\left(\frac{x}{x+1}\right)
+\frac{2}{{4\,k\,M}+4\,\ell+2}\left(\frac{x}{x+1}\right)^2
+\frac{2}{{4\,k\,M}+4\,\ell+3}\left(\frac{x}{x+1}\right)^3
\right)
\Biggr).
\end{align*}
\end{theorem}
$\\ $
The Proof is not difficult. Only a little Algebra and Classical Analysis are required.
\section{Proof of the theorem}
Let $r\in\mathbb{R}\backslash\left\{1\right\}$. First we state the identity
\[
\frac{1}{1-r}=\frac{1}{1-r}\cdot\frac{1-r^M}{1-r^M}=\frac{1}{1-r^M}\cdot\frac{1-r^M}{1-r}\;.
\]
From now on we assume $\left|r\right|<1$. Thus, the fraction $\frac{1}{1-{r}^{M}}$
can be expressed as an infinite geometric series \cite{AS}
\[
\frac {1}{1-{r}^{M}}
=\sum_{k=0}^{\infty}{\left(r^M\right)}^k=\sum_{k=0}^{\infty}{r}^{\,k\,M}
\]
The second Term $\frac{1-r^M}{1-r}$ can also be expressed as a finite geometric series \cite{SO}
\[
\frac{1-r^M}{1-r}=\sum_{l=0}^{M-1}{r}^{l}\;.
\]
Plugging in gives
\[
\frac{1}{1-r}=\frac{1}{1-r^M}\cdot\frac{1-r^M}{1-r}
=\left(\sum_{k=0}^{\infty}{r}^{k\,M}\right)\left(\sum_{\ell=0}^{M-1}{r}^{\ell}\right)
\]
We obtain a product of two sums which converges for $\left|r\right|<1$
\[
\frac{1}{1-r}
=\left(\sum_{k=0}^{\infty}{r}^{k\,M}\right)\left(\sum_{\ell=0}^{M-1}{r}^{\ell}\right)
\]
The right-hand side can easily represented as
\[
\frac{1}{1-r}=\left(\sum_{k=0}^{\infty}{r}^{k\,M}\right)\left(\sum_{\ell=0}^{M-1}{r}^{\ell}\right)
=\sum_{k=0}^{\infty}\left(\sum_{l=0}^{M-1}{r}^{k\,M+\ell}\right)
\;.
\]
Next, substituting
\[r=-4\,s^4\,.\]
This changes the interval of convergence $\left|r\right|<1$ to a new interval of convergence
\[\left|s\right|<\frac{1}{2}\sqrt{2}\,.\]
Plugging in $r=-4\,s^4$ in the formula above gives
\[
\frac{1}{1+s^4}=\sum_{k=0}^{\infty}\left(\sum_{\ell=0}^{M-1}\left(-4\,s^4\right)^{k\,M+\ell}\right).
\]
Expanding and rearranging the terms yields
\[
\frac{1}{1+s^4}
=\sum_{k=0}^{\infty}\left(\sum_{\ell=0}^{M-1}\left(-4\right)^{k\,M+\ell}s^{4\,k\,M+4\,\ell}\right).
\]
Using the statement, which can be easily proven
\[
1+4\,s^4=\left(1-2\,s+2\;s^2\right)\left(1+2\,s+2\;s^2\right)
\]
leads to
\begin{align*}
\frac{1}{\left(1-2\,s+2\;s^2\right)\left(1+2\,s+2\;s^2\right)}
=\sum_{k=0}^{\infty}\left(\sum_{\ell=0}^{M-1}\left(-4\right)^{k\,M+\ell}s^{4\,k\,M+4\,\ell}\right).
\end{align*}
Multiplying by
\[
\left(1+2\,s+2\;s^2\right)
\]
factoring out and reordering the terms results in
\begin{align*}
\frac{1}{1-2\,s+2\;s^2}
&=\left(1+2\,s+2\;s^2\right)
\sum_{k=0}^{\infty}\left(\sum_{\ell=0}^{M-1}\left(-4\right)^{k\,M+\ell}s^{4\,k\,M+4\,\ell}\right)
\\
&=\sum_{k=0}^{\infty}\left(\sum_{\ell=0}^{M-1}\left(-4\right)^{k\,M+\ell}
\left(s^{4\,k\,M+4\,\ell}+2\,s^{4\,k\,M+4\,\ell+1}+2\,s^{4\,k\,M+4\,\ell+2}\right)\right).
\end{align*}
Next step is integrating from $s=0$ to $s=t$. That gives
\begin{align*}
\int\limits_{0}^{t}\frac{1}{1-2\,s+2\;s^2}\,{\rm d}s
=\int\limits_{0}^{t}\sum_{k=0}^{\infty}
\left(\sum_{\ell=0}^{M-1}\left(-4\right)^{k\,M+\ell}
\left(
   s^{4\,k\,M+4\,\ell}
+2\,s^{4\,k\,M+4\,\ell+1}
+2\,s^{4\,k\,M+4\,\ell+2}
\right)\right)\,{\rm d}s
\end{align*}
We shall first turn our attention to the left-hand side. We use the antiderivative \cite{RC}
\[
\int\frac{1}{1-2\,s+2\;s^2}\,{\rm d}s=\arctan\left({2\,s-1}\right)+C\;
\]
and thus obtain
\begin{align*}
\int\limits_{0}^{t}\!\frac {1}{1-2\,s+2\;s^2}\,{\rm d}s
&=\Bigg[\,\arctan\left({2\,s-1}\right)\,\Bigg]_0^t
\\
&=\arctan\left({2\,t-1}\right)-\arctan\left(-1\right)
\\
&=\arctan\left({2\,t-1}\right)+\arctan\left(1\right)
\\
&=\arctan\left({2\,t-1}\right)+\frac{\pi}{4}\;.
\end{align*}
Let us now turn to the right-hand side
\begin{align*}
\int\limits_{0}^{t}
&\sum_{k=0}^{\infty}
\left(\sum_{\ell=0}^{M-1}\left(-4\right)^{k\,M+\ell}
\left(
    s^{4\,k\,M+4\,\ell}
+2\,s^{4\,k\,M+4\,\ell+1}
+2\,s^{4\,k\,M+4\,\ell+2}
\right)\right)\,{\rm d}s
\end{align*}
The absolute convergence allows the interchange of the sum and the integral \cite{KK}. That gives
\begin{align*}
\int\limits_{0}^{t}
&\sum_{k=0}^{\infty}
\left(\sum_{\ell=0}^{M-1}\left(-4\right)^{k\,M+\ell}
\left(
    s^{4\,k\,M+4\,\ell}
+2\,s^{4\,k\,M+4\,\ell+1}
+2\,s^{4\,k\,M+4\,\ell+2}
\right)\right)\,{\rm d}s
\\
=&\sum_{k=0}^{\infty}
\left(\sum_{\ell=0}^{M-1}\left(-4\right)^{k\,M+\ell}
\int\limits_{0}^{t}
\left(
    s^{4\,k\,M+4\,\ell}
+2\,s^{4\,k\,M+4\,\ell+1}
+2\,s^{4\,k\,M+4\,\ell+2}
\right)\,{\rm d}s\right)
\\
=&\sum_{k=0}^{\infty}
\left(\sum_{\ell=0}^{M-1}\left(-4\right)^{k\,M+\ell}
\left(
 \frac{t^{4\,k\,M+4\,\ell+1}}{4\,k\,M+4\,\ell+1}
+\frac{2\,t^{4\,k\,M+4\,\ell+2}}{4\,k\,M+4\,\ell+2}
+\frac{2\,t^{4\,k\,M+4\,\ell+3}}{4\,k\,M+4\,\ell+3}
\right)\right)
\end{align*}
Applying a little algebra, we get finally
\begin{align*}
\int\limits_{0}^{t}
&\sum_{k=0}^{\infty}
\left(\sum_{\ell=0}^{M-1}\left(-4\right)^{k\,M+\ell}
\left(
    s^{4\,k\,M+4\,\ell}
+2\,s^{4\,k\,M+4\,\ell+1}
+2\,s^{4\,k\,M+4\,\ell+2}
\right)\right)\,{\rm d}s
\\
=&\sum_{k=0}^{\infty}
\left(-4\,t^4\right)^{k\,M}\,
\left(\sum_{\ell=0}^{M-1}
\left(-4\,t^4\right)^{\ell}
\left(
    \frac{t}{k\,M+4\,\ell+1}
+\frac{2\;{t}^{2}}{k\,M+4\,\ell+2}
+\frac{2\;{t}^{3}}{k\,M+4\,\ell+3}
\right)\right)
\end{align*}
Up to now we have:
\begin{align*}
a)&\int\limits_{0}^{t}\frac{1}{1-2\,s+2\;s^2}\,{\rm d}s=\arctan\left({2\,t-1}\right)+\frac{\pi}{4}\;;
\\
b)&
\int\limits_{0}^{t}
\sum_{k=0}^{\infty}
\left(\sum_{\ell=0}^{M-1}\left(-4\right)^{k\,M+\ell}
\left(
    s^{4\,k\,M+4\,\ell}
+2\,s^{4\,k\,M+4\,\ell+1}
+2\,s^{4\,k\,M+4\,\ell+2}
\right)\right)\,{\rm d}s
\\
&\,=\sum_{k=0}^{\infty}
\left(-4\,t^4\right)^{k\,M}\,
\left(\sum_{\ell=0}^{M-1}
\left(-4\,t^4\right)^{\ell}
\left(
    \frac{t}{k\,M+4\,\ell+1}
+\frac{2\;{t}^{2}}{k\,M+4\,\ell+2}
+\frac{2\;{t}^{3}}{k\,M+4\,\ell+3}
\right)\right).
\end{align*}
Putting all the pieces together results in
\begin{align*}
&\arctan\left({2\,t-1}\right)+\frac{\pi}{4}
\\
&=\sum_{k=0}^{\infty}
\left(\left(-4\,t^4\right)^{M}\right)^k
\left(\sum_{\ell=0}^{M-1}\left(-4\,t^4\right)^{\ell}
\left(
\frac{t}{{4\,k\,M}+4\,\ell+1}
+\frac{2\;t^{2}}{{4\,k\,M}+4\,\ell+2}
+\frac{2\;t^{3}}{{4\,k\,M}+4\,\ell+3}
\right)\right)
\end{align*}
$\\ $
In this representations, we have $\frac{\pi}{4}$ on the left-hand side. It is therefore not particularly suitable for practical computation and is of more of an academic interest. However, this can easily be changed. Let us remind $\frac{\pi}{4}=\arctan\left(1\right)$ and the addition theorem for the real arctangent \cite{SO}
\[
\arctan\left(\alpha\right)+\arctan\left(\beta\right)
=\arctan\left(\frac{\alpha+\beta}{1-\alpha\,\beta}\right),\;
\alpha\,\beta < 1\;.
\]
Let us put
\begin{align*}
\alpha=2\,t-1\,;\quad\beta=1\,.
\end{align*}
That yields
\begin{align*}
\arctan\left({2\,t-1}\right)+\arctan\left(1\right)
&=\arctan\left(\frac{2\,t-1+1}{1-\bigl[\left(2\,t-1\right)\cdot 1\bigr]}\right)
\\
&=\arctan\left(\frac{2\,t}{2-2\,t}\right)
\\
&=\arctan\left(\frac{t}{1-t}\right)\,.
\end{align*}
The condition $\alpha\,\beta < 1$ implies that $t < 1$. However, since we have assumed that $\left|t\right|<\frac{1}{2}\sqrt{2}$, this condition is automatically fulfilled. We obtain a representation that is already rather useful
\begin{align*}
&\arctan\left(\frac{t}{1-t}\right)
\\
&=\sum_{k=0}^{\infty}
\left(\left(-4\,t^4\right)^{M}\right)^k
\left(\sum_{\ell=0}^{M-1}\left(-4\,t^4\right)^{\ell}
\left(
\frac{t}{{4\,k\,M}+4\,\ell+1}
+\frac{2\;t^{2}}{{4\,k\,M}+4\,\ell+2}
+\frac{2\;t^{3}}{{4\,k\,M}+4\,\ell+3}
\right)\right).
\end{align*}
It would, of course, be desirable to be able to calculate the value of $\arctan\left(x\right)$ directly for a given $x$. This can be done quite easily. We substitute
$t=\frac{x}{x+1}$ into $\arctan\left(\frac{t}{1-t}\right)$ and obtain on the left-hand side
\[
\arctan\left(\frac{t}{1-t}\right)
=\arctan\left(\frac{\frac{x}{x+1}}{1-\frac{x}{x+1}}\right)
=\arctan\left(\frac{\frac{x}{x+1}}{\frac{x+1}{x+1}-\frac{x}{x+1}}\right)
=\arctan\left(\frac{\frac{x}{x+1}}{\frac{1}{x+1}}\right)
=\arctan\left(x\right)\;.
\]
The left-hand side result is as desired. Substituting $t=\frac{x}{x+1}$ into the right-hand side and applying some algebra, we get
\begin{align*}
&\arctan\left(x\right)=\sum_{k=0}^{\infty}
\left(\left(\frac{-4\,x^4}{\left(x+1\right)^4}\right)^M\right)^{k}
\left(\sum_{\ell=0}^{M-1}\left(\frac{-4\,x^4}{\left(x+1\right)^4}\right)^{\ell}
\right.\times
\\
&\Biggl.
\left(
 \frac{1}{{4\,k\,M}+4\,\ell+1}\left(\frac{x}{x+1}\right)
+\frac{2}{{4\,k\,M}+4\,\ell+2}\left(\frac{x}{x+1}\right)^2
+\frac{2}{{4\,k\,M}+4\,\ell+3}\left(\frac{x}{x+1}\right)^3
\right)
\Biggr).
\end{align*}
Naturally, the substitution changes the interval of convergence.
The condition for convergence
\[\left|t\right|<\frac{1}{2}\sqrt{2}\]
becomes, following the substitution $t=\frac{x}{x+1}$,
\[
\left|\frac{x}{x+1}\right|< \frac{1}{2}\sqrt{2}\;.
\]
This inequality is easy to solve and yields the new interval of convergence, now for $x$
\[1-\sqrt{2}< x< 1+\sqrt{2}\;.\]
This completes the proof of the theorem.
\[\]
The figure below shows $\arctan\left(x\right)$ and some approximate curves.$\\ $
The interval of convergence
\[1-\frac{1}{2}\,\sqrt{2}< x<\frac{1}{2}\,\sqrt{2}\]
is also shown in the figure.
\[\]
\begin{center}
\fbox{\includegraphics[width=0.95\linewidth]{"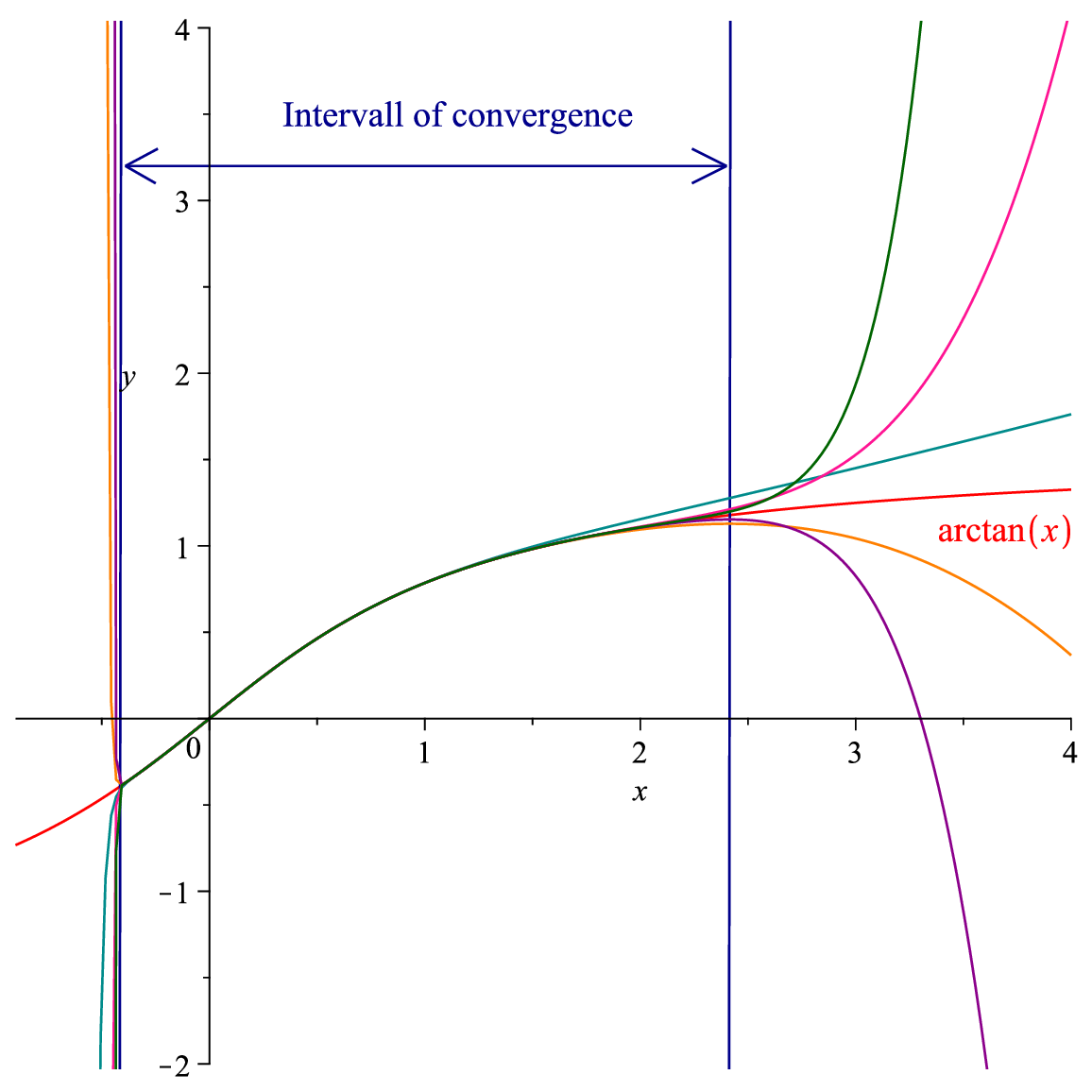"}}
\end{center}
$\\ $
\centerline{\textbf{\underline{figure 1\,:}} Approximation of $\arctan\left(x\right)$ and some series expansions.}
\section{\texorpdfstring{BBP Formulas for $\arctan\left(\frac{m}{n}\right)$}{}}
We now come to a very interesting application. We apply the series to $\arctan\left(\frac{m}{n}\right)$ wehre $m\in\mathbb{Z}$, $n\in\mathbb{N}$ and the statement $1-\sqrt{2}<\left|\frac{m}{n}\right|<1+\sqrt{2}$ holds true.
Using
\begin{align*}
&\arctan\left(x\right)
=\sum_{k=0}^{\infty }
\left(\left(\frac{-4\,x^4}{\left(x+1\right)^4}\right)^M\right)^k
\left(\sum_{\ell=0}^{M-1}\left(\frac{-4\,x^4}{\left(x+1\right)^4}\right)^{\ell}
\right.\times
\\
&\Biggl.
\left(
 \frac{1}{{4\,k\,M}+4\,\ell+1}\left(\frac{x}{x+1}\right)
+\frac{2}{{4\,k\,M}+4\,\ell+2}\left(\frac{x}{x+1}\right)^2
+\frac{2}{{4\,k\,M}+4\,\ell+3}\left(\frac{x}{x+1}\right)^3
\right)
\Biggr)
\end{align*}
and plugging in $x=\frac{m}{n}$ leads to
\begin{align*}
&\arctan \left(\frac{m}{n}\right)
=\sum_{k=0}^{\infty}\left(\left(\frac{-4\,m^4}{\left(m+n\right)^4}\right)^{M}\right)^k
\;\left(\sum_{\ell=0}^{M-1}\left(\frac{-4\,m^4}{\left(m+n\right)^4}\right)^{\ell}\right.\times
\\
&\Bigg.
\left(
 \frac{m}{\left(4\,k\,M+4\,\ell+1\right)\left(m+n\right)}
+\frac{2\,m^2}{\left(4\,k\,M+4\,\ell+2\right)\left(m+n\right)^2}
+\frac{2\,m^3}{\left(4\,k\,M+4\,\ell+3\right)\left(m+n\right)^3}
\right)
\Bigg)
\end{align*}
We identify the series as series of the BBP’ type with base
\[b=\left(\frac{-4\,m^4}{\left(m+n\right)^4}\right)^{M}.\]
Every $M\in\mathbb N$ gives rise to a new series. BBP formulas are of interest because, for fixed $b$, the $k^{th}$ $b$-ary digit of a number with a BBP formula can be found without knowing prior digits.
\subsection{\texorpdfstring{Machin-Type BBP Formulas}{}}
Machin’s formulas are particular identities for $\frac{\pi}{4}$ as a sum of integer multiples of arctangent having rational arguments
\[
\frac{\pi}{4}=\sum_{k=1}^{L}a_{\,k}\arctan\frac{1}{b_{\,k}}\;,
\]
where $a_{\,k}\in\mathbb Z$ and $b_{\,k}\in\mathbb N$.
Here, the argument takes the value $x=\frac{1}{n}$. That means $m=1$. Plugging in $m=1$ into the formula above leads to
\begin{align*}
&\arctan \left(\frac{1}{n}\right)
=\sum_{k=0}^{\infty}\left(\left(\frac{-4}{\left(n+1\right)^4}\right)^{M}\right)^k
\;\left(\sum_{\ell=0}^{M-1}\left(\frac{-4}{\left(n+1\right)^4}\right)^{\ell}\right.\times
\\
&\Bigg.
\left(
 \frac{1}{\left(4\,k\,M+4\,\ell+1\right)\left(n+1\right)}
+\frac{2}{\left(4\,k\,M+4\,\ell+2\right)\left(n+1\right)^2}
+\frac{2}{\left(4\,k\,M+4\,\ell+3\right)\left(n+1\right)^3}
\right)
\Bigg)
\end{align*}
Now we have a simpler base
\[b=\left(\frac{-4}{\left(n+1\right)^4}\right)^{M}.\]
\underline{Note}
$\\ $
The special case $n=2^m-1\,,\;m\in\mathbb{N}$ leads to a binary base
\[
b=\left(\frac{-4}{\left(2^{\,m}\right)^4}\right)^{M}
=\left(\frac{-4}{2^{\,4\,m}}\right)^{M}
=\left(\frac{-1}{2^{\,4\,m-2}}\right)^{M}.
\]

\subsubsection{1-Term formula}
There is a single formula with one term
\[
\frac{\pi}{4}=\arctan\left(1\right)
\]
So $n=1$, we obtain
\begin{align*}
&\frac{\pi}{4}
=\left(\sum_{k=0}^{\infty}\left(\left(-\frac{1}{4}\right)^M\right)^k\times\right.
\\
&\left(\sum_{\ell=0}^{M-1}\left(-\frac{1}{4}\right)^\ell
\left(
 \frac{1}{2\left({4\,k\,M}+4\,\ell+1\right)}
+\frac{1}{2\left({4\,k\,M}+4\,\ell+2\right)}
+\frac{1}{4\left({4\,k\,M}+4\,\ell+3\right)}
\right)\right)
\end{align*}
$\\ $
In this special case we have base $\left(-\frac{1}{4}\right)^M$ and therefore binary digit extraction formulas for $\frac{\pi}{4}$.
$\\ $
Let’s take a look at the formulas for some $M$

\[\underline{ M=1}\]
\[\frac{\pi}{4}=\sum_{k=0}^{\infty}
\left(-\frac{1}{4}\right)^k
\left(
 \frac{1}{2\,\left(4\,k+1\right)}
+\frac{1}{2\,\left(4\,k+2\right)}
+\frac{1}{4\,\left(4\,k+3\right)}
\right)
\]
$\\ $
This formula was published in 1997 by Victor Adamchik and Stan Wagon \cite{AW}. It was already known to Euler \cite{LE}, albeit in a slightly different notation (Euler used planetary symbols for the variables).
\[\underline{ M=2}\]
\[\frac{\pi}{4}=\sum_{k=0}^{\infty}\left(\frac{1}{16}\right)^k\times\]
\[
\left(
 \frac{1}{2\,\left(8\,k+1\right)}
+\frac{1}{2\,\left(8\,k+2\right)}
+\frac{1}{4\,\left(8\,k+3\right)}
\right.
\]
\[
\left.
-\frac{1}{8\,\left(8\,k+5\right)}
-\frac{1}{8\,\left(8\,k+6\right)}
-\frac{1}{16\,\left(8\,k+7\right)}
\right)
\]

\[\underline{ M=3}\]
\[\frac{\pi}{4}=\sum_{k=0}^{\infty}\left(-\frac{1}{64}\right)^k\times\]
\[
\left(
 \frac{1}{2\,\left(12\,k+1\right)}
+\frac{1}{2\,\left(12\,k+2\right)}
+\frac{1}{4\,\left(12\,k+3\right)}
\right.
\]
\[
-\frac{1}{8\,\left(12\,k+5\right)}
-\frac{1}{8\,\left(12\,k+6\right)}
-\frac{1}{16\,\left(12\,k+7\right)}
\]
\[
\left.
+\frac{1}{32\,\left(12\,k+ 9\right)}
+\frac{1}{32\,\left(12\,k+10\right)}
+\frac{1}{64\,\left(12\,k+11\right)}
\right)
\]

\[\underline{ M=4}\]
\[\frac{\pi}{4}=\sum_{k=0}^{\infty}\left(\frac{1}{256}\right)^k\times
\]
\[
\left(
 \frac{1}{2\,\left(16\,k+1\right)}
+\frac{1}{2\,\left(16\,k+2\right)}
+\frac{1}{4\,\left(16\,k+3\right)}
\right.
\]
\[
-\frac{1}{8\,\left(16\,k+5\right)}
-\frac{1}{8\,\left(16\,k+6\right)}
-\frac{1}{16\,\left(16\,k+7\right)}
\]
\[
+\frac{1}{32\,\left(16\,k+ 9\right)}
+\frac{1}{32\,\left(16\,k+10\right)}
+\frac{1}{64\,\left(16\,k+11\right)}
\]
\[
\left.
-\frac{1}{128\,\left(16\,k+13\right)}
-\frac{1}{128\,\left(16\,k+14\right)}
-\frac{1}{256\,\left(16\,k+15\right)}
\right)
\]
\[\vdots\]
%
%
\subsection{2-Term formulas}
According to Størmer \cite{CS}, there are exactly four formulas with two terms
\begin{align*}
&\frac{\pi}{4}=\arctan\frac{1}{2}+\arctan\frac{1}{3}
\\
&\frac{\pi}{4}=2\arctan\frac{1}{2}-\arctan\frac{1}{7}
\\
&\frac{\pi}{4}=2\arctan\frac{1}{3}+\arctan\frac{1}{7}
\\
&\frac{\pi}{4}=4\arctan\frac{1}{5}-\arctan\frac{1}{239}
\end{align*}
We start with
\[\frac{\pi}{4}
=\arctan\left(\frac{1}{2}\right)+\arctan\left(\frac{1}{3}\right)
\]
For $\arctan\left(\frac{1}{2}\right)$, that is, $n=2$, we have
\begin{align*}
&\arctan\left(\frac{1}{2}\right)
=\sum_{k=0}^{\infty}\left(\left(-{\frac{4}{81}}\right)^{M}\right)^{k}\times
\\
&\Biggl(\sum_{\ell=0}^{M-1}\left(-{\frac{4}{81}}\right)^{\ell}
\left(
 \frac{1}{3\left(4\,k\,M+4\,\ell+1\right)}
+\frac{2}{9\left(4\,k\,M+4\,\ell+2\right)}
+\frac{2}{27\left(4\,k\,M+4\,\ell+3\right)}
\right)
\Biggr)
\end{align*}

Evaluating $\arctan\left(\frac{1}{3}\right)$ yields
\begin{align*}
&\arctan\left(\frac{1}{3}\right)
=\sum_{k=0}^{\infty}\left(\left(-{\frac{1}{64}}\right)^{M}\right)^{k}\times
\\
&\Biggl(\sum_{\ell=0}^{M-1}\left(-{\frac{1}{64}}\right)^{\ell}
\left(
 \frac{1}{ 4\left(4\,k\,M+4\,\ell+1\right)}
+\frac{1}{ 8\left(4\,k\,M+4\,\ell+2\right)}
+\frac{1}{32\left(4\,k\,M+4\,\ell+3\right)}
\right)
\Biggr)
\end{align*}
$ \\ $
The next two formulas are
\[\frac{\pi}{4}=2\arctan\left(\frac{1}{2}\right)-\arctan\left(\frac{1}{7}\right);\]
\[\frac{\pi}{4}=2\arctan\left(\frac{1}{3}\right)+\arctan\left(\frac{1}{7}\right).\]
We already have the formula for $\arctan\left(\frac{1}{2}\right)$ and $\arctan\left(\frac{1}{3}\right)$, so we only need them for $\arctan\left(\frac{1}{7}\right)$
\begin{align*}
&\arctan\left(\frac{1}{7}\right)
=\sum_{k=0}^{\infty}\left(\left(-\frac{1}{1024}\right)^{M}\right)^{k}\times
\\
&\Biggl(\sum_{\ell=0}^{M-1}\left(-\frac{1}{1024}\right)^{\ell}
\left(
 \frac{1}{  8\left(4\,k\,M+4\,\ell+1\right)}
+\frac{1}{ 32\left(4\,k\,M+4\,\ell+2\right)}
+\frac{1}{256\left(4\,k\,M+4\,\ell+3\right)}
\right)
\Biggr).
\end{align*}
$ \\ $
For $M=1$, $\arctan\left(\frac{1}{3}\right)$ and $\arctan\left(\frac{1}{7}\right)$ yield series that were already known to Euler \cite{LE}.
$ \\ $
Now we come to Machin’s famous formula
\[\frac{\pi}{4}
=4\arctan\left(\frac{1}{5}\right)-\arctan\left(\frac{1}{239}\right)
\]
For $\arctan\left(\frac{1}{5}\right)$, we get
\begin{align*}
&\arctan\left(\frac{1}{5}\right)=
\sum_{k=0}^{\infty}\left(\left(-\frac{1}{324}\right)^{M}\right)^{k}\times
\\
&\Biggl(\sum_{\ell=0}^{M-1}\left(-{\frac{1}{324}}\right)^{\ell}
\left(
 \frac{1}{ 6\left(4\,k\,M+4\,\ell+1\right)}
+\frac{1}{18\left(4\,k\,M+4\,\ell+2\right)}
+\frac{1}{108\left(4\,k\,M+4\,\ell+3\right)}
\right)
\Biggr)
\end{align*}

Finally $\arctan\left(\frac{1}{239}\right)$ yields
\begin{align*}
&\arctan\left(\frac{1}{239}\right)=\sum_{k=0}^{\infty}
\left(\left(-{\frac{1}{829440000}}\right)^{M}\right)^{k}\times
\Biggl(\sum_{\ell=0}^{M-1}\left(-{\frac{1}{829440000}}\right)^{\ell}
\\
&\Biggl.
\left(
 \frac{1}{240\left(4\,k\,M+4\,\ell+1\right)}
+\frac{1}{28800\left(4\,k\,M+4\,\ell+2\right)}
+\frac{1}{6912000\left(4\,k\,M+4\,\ell+3\right)}
\right)
\Biggr)
\end{align*}
\subsection{formulas with more then two terms}
The next step involves Machin-Type formulas with more than two terms. We shall give the most effective ones with three and four terms, leaving it to the reader to work out the series.
\[\frac{\pi}{4}
=12\arctan\left(\frac{1}{18}\right)
+8\arctan\left(\frac{1}{57}\right)
-5\arctan\left(\frac{1}{239}\right)
\]
\[\frac{\pi}{4}
=44\arctan\left(\frac{1}{57}\right)
+7\arctan\left(\frac{1}{239}\right)
-12\arctan\left(\frac{1}{682}\right)
+24\arctan\left(\frac{1}{12943}\right)
\] 
For further $\frac{\pi}{4}$ formulas of the Machin type, we refer the reader to the very detailed paper by Jörg Arndt \cite{JA}. The more terms are used, the smaller the argument becomes and the shorter the computation time. For example, if $\arctan\left(\frac{1}{n}\right)$ is to be computed with an already very small argument $\frac{1}{n}$, then using the method presented here results in an even smaller argument
\[\left(\frac{-4}{(n+1)^4}\right)^{M}.\]

{}

\newpage
\selectlanguage{ngerman}
\begin{center}
\LARGE{{Eine rationale Reihendarstellung für den Arkustangens und ein Verfahren zur Konstruktion von BBP Reihen für Formeln vom Machin'schen Typ
}}
\[\]
\[\]
\[\]
\[\]
\[\]
\[\]
\large
{Alois Schiessl}
\footnote[1]{\,\;E-Mail: \texttt{aloisschiessl@posteo.de}}
\centerline{ }
\centerline{\today}
\end{center}
\[\]
\[\]
\[\]
\[\]
\begin{abstract}
In dieser Veröffentlichung geben wir für die reelle $\arctan\left(x\right)$ Funktion eine Reihendarstellung an, die aus einer Summe aus gebrochen rationalen Funktionen besteht, und einen frei wählbaren Parameter $M\in \mathbb{N}$ enthält. Jede natürliche Zahl $M$ führt zu einer neuen Reihe. Damit können effiziente BBP-Reihen für Formeln vom Machin'schen Typ konstruiert werden. Speziell für $x=1$ ergeben sich Reihen für die binäre Bit-Extraktion.
\end{abstract}

\setcounter{section}{0}
\section{Einführung und Hauptergebnis}
In dieser Veröffentlichung geben wir für die reelle $\arctan\left(x\right)$ Funktion eine Darstellung an, die aus einer Summe aus gebrochen rationalen Funktionen besteht, und die einen frei wählbaren Parameter $M\in \mathbb{N}$ enthält. Jede natürliche Zahl $M$ führt zu einer neuen Reihe. Damit können effiziente BBP-Reihen für Formeln vom Machin'schen Typ konstruiert werden. Speziell für $x=1$ ergeben sich Reihen für die binäre Bit-Extraktion.
\newpage
$ \\ $
Wir formulieren unsere Aussage
\begin{theorem}
$ \\ $
Es sei $1-\sqrt{2}< x < 1+\sqrt{2}$ und $M\in\mathbb N$. Dann gilt die Reihendarstellung:
\begin{align*}
&\arctan\left(x\right)=\sum_{k=0}^{\infty}
\left(\left(\frac{-4\,x^4}{\left(x+1\right)^4}\right)^M\right)^{k}
\left(\sum_{\ell=0}^{M-1}\left(\frac{-4\,x^4}{\left(x+1\right)^4}\right)^{\ell}
\right.\times
\\
&\Biggl.
\left(
 \frac{1}{{4\,k\,M}+4\,\ell+1}\left(\frac{x}{x+1}\right)
+\frac{2}{{4\,k\,M}+4\,\ell+2}\left(\frac{x}{x+1}\right)^2
+\frac{2}{{4\,k\,M}+4\,\ell+3}\left(\frac{x}{x+1}\right)^3
\right)
\Biggr).
\end{align*}
\end{theorem}
Der Beweis hierfür ist nicht schwierig. Es werden nur elementare Algebra und klassische Analysis benötigt.

\section{Beweis Theorem}
Es $r\in \mathbb{R}\backslash \left\{1\right\}$. Dann gilt die Identität
\[
\frac{1}{1-r}=\frac{1}{1-r}\cdot\frac{1-r^M}{1-r^M}=\frac{1}{1-r^M}\cdot\frac{1-r^M}{1-r}\;.
\]
Von nun an werde $\left|r\right|<1$ vorausgesetzt. Somit ist der Bruch $\frac{1}{1-{r}^{M}}$
durch eine unendliche geometrische Reihe darstellbar \cite{AS}
\[
\frac {1}{1-{r}^{M}}
=\sum_{k=0}^{\infty}{\left(r^M\right)}^k=\sum_{k=0}^{\infty}{r}^{\,k\,M}
\]
Wir erkennen den Term $\frac{1-r^M}{1-r}$ als Summenformel für die endliche geometrische Reihe
\[
\frac{1-r^M}{1-r}=\sum_{l=0}^{M-1}{r}^{l}\;.
\]
Wir erhalten dadurch ein Produkt aus zwei Summen mit gemeinsamen Konvergenzbereich $\left|r\right|<1$
\[
\frac{1}{1-r}
=\left(\sum_{k=0}^{\infty}{r}^{k\,M}\right)\left(\sum_{\ell=0}^{M-1}{r}^{\ell}\right)
\]
Die Summe
\[\sum_{l=0}^{M-1}{r}^{l}\]
hat für festes $r$ und $M$ einen endlichen Wert und kann somit in die unendliche Summe gebracht werden. Das Produkt der Summen lässt sich damit auf einfache Weise wie folgt darstellen:
\[
\frac{1}{1-r}=\left(\sum_{k=0}^{\infty}{r}^{k\,M}\right)\left(\sum_{\ell=0}^{M-1}{r}^{\ell}\right)
=\sum_{k=0}^{\infty}\left(\sum_{l=0}^{M-1}{r}^{k\,M+\ell}\right)
\;.
\]
Als nächstes substituieren wir
$r=-4\,s^4\;$.
Die Konvergenzbedingung $\left|r\right|<1$ führt zur neuen Konvergenzbedingung 
\[\left|s\right|<\frac{1}{2}\sqrt{2}\]
Die Substitution $r=-4\,s^4$ ergibt dann
\[
\frac{1}{1+s^4}=\sum_{k=0}^{\infty}\left(\sum_{\ell=0}^{M-1}\left(-4\,s^4\right)^{k\,M+\ell}\right)
\]
Das lässt sich weiter umformen zu
\[
\frac{1}{1+s^4}
=\sum_{k=0}^{\infty}\left(\sum_{\ell=0}^{M-1}\left(-4\right)^{k\,M+\ell}s^{4\,k\,M+4\,\ell}\right)
\]
Nun gilt folgende Zerlegung:
\[
1+4\,s^4=\left(1-2\,s+2\;s^2\right)\left(1+2\,s+2\;s^2\right)
\]
Das lässt sich durch direktes Ausrechnen bestätigen. Wir können also auch schreiben
\begin{align*}
\frac{1}{\left(1-2\,s+2\;s^2\right)\left(1+2\,s+2\;s^2\right)}
=\sum_{k=0}^{\infty}\left(\sum_{\ell=0}^{M-1}\left(-4\right)^{k\,M+\ell}s^{4\,k\,M+4\,\ell}\right)
\end{align*}
Als nächstes multiplizieren wir mit
\[
\left(1+2\,s+2\;s^2\right)\;.
\]
Das ergibt dann
\begin{align*}
\frac{1}{1-2\,s+2\;s^2}
&=\left(1+2\,s+2\;s^2\right)
\sum_{k=0}^{\infty}\left(\sum_{\ell=0}^{M-1}\left(-4\right)^{k\,M+\ell}s^{4\,k\,M+4\,\ell}\right)
\end{align*}
Das Polynom $1+2\,s+2\;s^2$ ist von endlicher Größe. Es kann wie eine Konstante in die Summe gezogen werden. Wir erhalten somit
\begin{align*}
\frac{1}{1-2\,s+2\;s^2}
=\sum_{k=0}^{\infty}\left(\sum_{\ell=0}^{M-1}\left(-4\right)^{k\,M+\ell}
\left(s^{4\,k\,M+4\,\ell}+2\,s^{4\,k\,M+4\,\ell+1}+2\,s^{4\,k\,M+4\,\ell+2}\right)\right).
\end{align*}
Im nächsten Schritt integrieren wir zwischen den Grenzen $s=0$ und $s=t$ und erhalten
\begin{align*}
\int\limits_{0}^{t}\frac{1}{1-2\,s+2\;s^2}\,{\rm d}s
=\int\limits_{0}^{t}\sum_{k=0}^{\infty}
\left(\sum_{\ell=0}^{M-1}\left(-4\right)^{k\,M+\ell}
\left(
   s^{4\,k\,M+4\,\ell}
+2\,s^{4\,k\,M+4\,\ell+1}
+2\,s^{4\,k\,M+4\,\ell+2}
\right)\right)\,{\rm d}s
\end{align*}
Wir wenden uns zuerst der linken Seite zu. Die Verwendung der Stammfunktion \cite{WG}
\[
\int\frac{1}{1-2\,s+2\;s^2}\,{\rm d}s=\arctan\left({2\,s-1}\right)+C\;,
\]
die sich leicht durch Differentiation beweisen lässt, führt zu
\begin{align*}
\int\limits_{0}^{t}\!\frac {1}{1-2\,s+2\;s^2}\,{\rm d}s
&=\Bigg[\,\arctan\left({2\,s-1}\right)\,\Bigg]_0^t
\\
&=\arctan\left({2\,t-1}\right)-\arctan\left(-1\right)
\\
&=\arctan\left({2\,t-1}\right)+\arctan\left(1\right)
\\
&=\arctan\left({2\,t-1}\right)+\frac{\pi}{4}\;.
\end{align*}
Nun wenden wir uns der rechten Seite zu. Wir erhalten das Integral einer unendlichen Summe
\begin{align*}
\int\limits_{0}^{t}
&\sum_{k=0}^{\infty}
\left(\sum_{\ell=0}^{M-1}\left(-4\right)^{k\,M+\ell}
\left(
    s^{4\,k\,M+4\,\ell}
+2\,s^{4\,k\,M+4\,\ell+1}
+2\,s^{4\,k\,M+4\,\ell+2}
\right)\right)\,{\rm d}s
\end{align*}
Für $\left|s\right|<\frac{1}{2}\,\sqrt{2}\;$ liegt absolute Konvergenz vor, somit sind Integral und Summe vertauschbar \cite{HJ}
\begin{align*}
\int\limits_{0}^{t}
&\sum_{k=0}^{\infty}
\left(\sum_{\ell=0}^{M-1}\left(-4\right)^{k\,M+\ell}
\left(
    s^{4\,k\,M+4\,\ell}
+2\,s^{4\,k\,M+4\,\ell+1}
+2\,s^{4\,k\,M+4\,\ell+2}
\right)\right)\,{\rm d}s
\\
=&\sum_{k=0}^{\infty}
\left(\sum_{\ell=0}^{M-1}\left(-4\right)^{k\,M+\ell}
\int\limits_{0}^{t}
\left(
    s^{4\,k\,M+4\,\ell}
+2\,s^{4\,k\,M+4\,\ell+1}
+2\,s^{4\,k\,M+4\,\ell+2}
\right)\,{\rm d}s\right)
\\
=&\sum_{k=0}^{\infty}
\left(\sum_{\ell=0}^{M-1}\left(-4\right)^{k\,M+\ell}
\left(
 \frac{t^{4\,k\,M+4\,\ell+1}}{4\,k\,M+4\,\ell+1}
+\frac{2\,t^{4\,k\,M+4\,\ell+2}}{4\,k\,M+4\,\ell+2}
+\frac{2\,t^{4\,k\,M+4\,\ell+3}}{4\,k\,M+4\,\ell+3}
\right)\right)
\end{align*}
Bis jetzt haben wir folgende Ergebnisse:
\begin{align*}
a)&\int\limits_{0}^{t}\frac{1}{1-2\,s+2\;s^2}\,{\rm d}s=\arctan\left({2\,t-1}\right)+\frac{\pi}{4}\;;
\\
b)&\int\limits_{0}^{t}
\sum_{k=0}^{\infty}
\left(\sum_{\ell=0}^{M-1}\left(-4\right)^{k\,M+\ell}
\left(
    s^{4\,k\,M+4\,\ell}
+2\,s^{4\,k\,M+4\,\ell+1}
+2\,s^{4\,k\,M+4\,\ell+2}
\right)\right)\,{\rm d}s
\\
&=\sum_{k=0}^{\infty}
\left(\sum_{\ell=0}^{M-1}\left(-4\right)^{k\,M+\ell}
\left(
 \frac{t^{4\,k\,M+4\,\ell+1}}{4\,k\,M+4\,\ell+1}
+\frac{2\,t^{4\,k\,M+4\,\ell+2}}{4\,k\,M+4\,\ell+2}
+\frac{2\,t^{4\,k\,M+4\,\ell+3}}{4\,k\,M+4\,\ell+3}
\right)\right).
\end{align*}
Wir können den Term $\left(-4\,t^4\right)^{M}$ aus der inneren Summe herausziehen und vor die innere Summe stellen. Des Weiteren können wir $\left(-4\,t^4\right)$ in der inneren Summe ausklammern. Wir erhalten dann
\begin{align*}
\int\limits_{0}^{t}
&\sum_{k=0}^{\infty}
\left(\sum_{\ell=0}^{M-1}\left(-4\right)^{k\,M+\ell}
\left(
    s^{4\,k\,M+4\,\ell}
+2\,s^{4\,k\,M+4\,\ell+1}
+2\,s^{4\,k\,M+4\,\ell+2}
\right)\right)\,{\rm d}s
\\
=&\sum_{k=0}^{\infty}
\left(\left(-4\,t^4\right)^{M}\right)^k
\left(\sum_{\ell=0}^{M-1}\left(-4\,t^4\right)^{\ell}
\left(
\frac{t}{{4\,k\,M}+4\,\ell+1}
+\frac{2\;t^{2}}{{4\,k\,M}+4\,\ell+2}
+\frac{2\;t^{3}}{{4\,k\,M}+4\,\ell+3}
\right)\right)
\end{align*}
Fügen wir alle Teile zusammen, so ergibt sich eine Reihen-Darstellung für den Arkustangens mit Argument $\left(2\,t-1\right)$, gültig für $\left|t\right|<\frac{1}{2}\sqrt{2}$.
\begin{align*}
&\arctan\left({2\,t-1}\right)+\frac{\pi}{4}
\\
&=\sum_{k=0}^{\infty}
\left(\left(-4\,t^4\right)^{M}\right)^k
\left(\sum_{\ell=0}^{M-1}\left(-4\,t^4\right)^{\ell}
\left(
\frac{t}{{4\,k\,M}+4\,\ell+1}
+\frac{2\;t^{2}}{{4\,k\,M}+4\,\ell+2}
+\frac{2\;t^{3}}{{4\,k\,M}+4\,\ell+3}
\right)\right)
\end{align*}
In dieser Reihe tritt auf der linken Seite $\frac{\pi}{4}$ auf. Somit ist sie für praktische Berechnungen nicht so geeignet und eher vom akademischen Interesse. Das lässt sich aber leicht beseitigen. Wir erinnern wir uns an $\frac{\pi}{4}=\arctan\left(1\right)$ und das Additionstheorem für den reellen Arkustangens \cite{SO}
\[
\arctan\left(\alpha\right)+\arctan\left(\beta\right)
=\arctan\left(\frac{\alpha+\beta}{1-\alpha\,\beta}\right),\;
\alpha\,\beta < 1\;.
\]
Wir setzen:
$\alpha=2\,t-1\,;\quad\beta=1\,.$
$\\ $
$\\ $
Das ergibt dann:
\begin{align*}
\arctan\left({2\,t-1}\right)+\arctan\left(1\right)
&=\arctan\left(\frac{2\,t-1+1}{1-\bigl[\left(2\,t-1\right)\cdot 1\bigr]}\right)
\\
&=\arctan\left(\frac{2\,t}{2-2\,t}\right)
\\
&=\arctan\left(\frac{t}{1-t}\right)
\end{align*}
Die Bedingung $\alpha\beta<1$ führt zu $t<1$. Da wir aber $\left|t\right|<\frac{1}{2}\sqrt{2}$ vorausgesetzt haben, ist diese Forderung automatisch erfüllt. Wir erhalten eine schon eher brauchbare Darstellung
\begin{align*}
&\arctan\left(\frac{t}{1-t}\right)
\\
&=\sum_{k=0}^{\infty}
\left(\left(-4\,t^4\right)^{M}\right)^k
\left(\sum_{\ell=0}^{M-1}\left(-4\,t^4\right)^{\ell}
\left(
\frac{t}{{4\,k\,M}+4\,\ell+1}
+\frac{2\;t^{2}}{{4\,k\,M}+4\,\ell+2}
+\frac{2\;t^{3}}{{4\,k\,M}+4\,\ell+3}
\right)\right)
\end{align*}
Es wäre natürlich wünschenswert, wenn man zu einem vorgegebenen $x$ den Wert $\arctan\left(x\right)$ direkt berechnen könnte. Das ist ganz einfach möglich. Wir setzen 
$t=\frac{x}{x+1}$ in $\arctan\left(\frac{t}{1-t}\right)$ ein und erhalten auf der linken Seite
\[
\arctan\left(\frac{t}{1-t}\right)
=\arctan\left(\frac{\frac{x}{x+1}}{1-\frac{x}{x+1}}\right)
=\arctan\left(\frac{\frac{x}{x+1}}{\frac{x+1}{x+1}-\frac{x}{x+1}}\right)
=\arctan\left(\frac{\frac{x}{x+1}}{\frac{1}{x+1}}\right)
=\arctan\left(x\right)\;.
\]
Die linke Seite passt schon einmal. Einsetzen von $t=\frac{x}{x+1}$ in die rechte Seite der obigen Reihe ergibt nach einigen algebraischen Umformungen die im Theorem angegebene Reihendarstellung
\begin{align*}
&\arctan\left(x\right)=\sum_{k=0}^{\infty}
\left(\left(\frac{-4\,x^4}{\left(x+1\right)^4}\right)^M\right)^{k}
\left(\sum_{\ell=0}^{M-1}\left(\frac{-4\,x^4}{\left(x+1\right)^4}\right)^{\ell}
\right.\times
\\
&\Biggl.
\left(
 \frac{1}{{4\,k\,M}+4\,\ell+1}\left(\frac{x}{x+1}\right)
+\frac{2}{{4\,k\,M}+4\,\ell+2}\left(\frac{x}{x+1}\right)^2
+\frac{2}{{4\,k\,M}+4\,\ell+3}\left(\frac{x}{x+1}\right)^3
\right)
\Biggr).
\end{align*}
Das ist eine Monsterreihe, deren Summanden aus gebrochen rationalen Funktionen bestehen.
$\\ $
Natürlich ändert sich dadurch das Konvergenzverhalten.
Die Voraussetzung
\[\left|t\right|<\frac{1}{2}\sqrt{2}\]
wird wird durch die Substitution $t=\frac{x}{x+1}$ zu
\[
\left|\frac{x}{x+1}\right|< \frac{1}{2}\sqrt{2}\;.
\]
Diese Ungleichung ist leicht aufzulösen und ergibt die neue Konvergenzbedingung, jetzt für $x$:
\[1-\sqrt{2}< x< 1+\sqrt{2}\;.\]
Damit ist der Beweis für das Theorem abgeschlossen.
$\\ $
$\\ $
Die nachfolgende Abbildung zeigt $\arctan\left(x\right)$ und und einige Näherungskurven.
$\\ $
Das Konvergenz-Intervall
\[1-\frac{1}{2}\,\sqrt{2}< x< 1+\frac{1}{2}\,\sqrt{2}\]
ist ebenfalls mit abgebildet.
\[\]
\begin{center}
\fbox{\includegraphics[width=0.95\linewidth]{"M1M5.eps"}}
\end{center}
$\\ $
\centerline{\textbf{\underline{Abbildung 1\,:}} Approximation von $\arctan\left(x\right)$ durch Näherungssummen.}
\newpage
\section{\texorpdfstring{BBP Formeln für $\arctan\left(\frac{m}{n}\right)$}{}}
Jetzt kommen wir zu einer sehr interessanten Anwendung. Es sei $m\in\mathbb{Z}$, $n\in\mathbb{N}$. Wir bilden den Bruch $\frac{m}{n}$, und es sei die Konvergenzbedingung $1-\sqrt{2}<\left|\frac{m}{n}\right|<1+\sqrt{2}$ erfüllt. Wir setzen $x=\frac{m}{n}$ ein in die Reihe
\begin{align*}
&\arctan\left(x\right)
=\sum_{k=0}^{\infty }
\left(\left(\frac{-4\,x^4}{\left(x+1\right)^4}\right)^M\right)^k
\left(\sum_{\ell=0}^{M-1}\left(\frac{-4\,x^4}{\left(x+1\right)^4}\right)^{\ell}
\right.\times
\\
&\Biggl.
\left(
 \frac{1}{{4\,k\,M}+4\,\ell+1}\left(\frac{x}{x+1}\right)
+\frac{2}{{4\,k\,M}+4\,\ell+2}\left(\frac{x}{x+1}\right)^2
+\frac{2}{{4\,k\,M}+4\,\ell+3}\left(\frac{x}{x+1}\right)^3
\right)
\Biggr)
\end{align*}
Das ergibt dann nach einigen algebraischen Umformungen
\begin{align*}
&\arctan \left(\frac{m}{n}\right)
=\sum_{k=0}^{\infty}\left(\left(\frac{-4\,m^4}{\left(m+n\right)^4}\right)^{M}\right)^k
\;\left(\sum_{\ell=0}^{M-1}\left(\frac{-4\,m^4}{\left(m+n\right)^4}\right)^{\ell}\right.\times
\\
&\Bigg.
\left(
 \frac{m}{\left(4\,k\,M+4\,\ell+1\right)\left(m+n\right)}
+\frac{2\,m^2}{\left(4\,k\,M+4\,\ell+2\right)\left(m+n\right)^2}
+\frac{2\,m^3}{\left(4\,k\,M+4\,\ell+3\right)\left(m+n\right)^3}
\right)
\Bigg)
\end{align*}
Wir erkennen die Reihe als Reihe vom BBP’schen Typ mit der Basis
\[b=\left(\frac{-4\,m^4}{\left(m+n\right)^4}\right)^{M}.\]
Jedes $M\in\mathbb N$ erzeugt eine neue Reihe. BBP sind deshalb so interessant, weil man mit ihnen für eine feste Basis $b$ an einer beliebigen Stelle $k$ die sich dort befindliche Zahl berechnen kann, ohne dass man alle vorhergehenden Stellen kennen muss.
\subsection{BBP-Formeln für Reihen vom Machin'schen Typ}
Wir wenden die Reihe auf Formeln vom Machin' schen Typ an. Das sind Reihen für $\frac{\pi}{4}$ der speziellen Gestalt
\[
\frac{\pi}{4}=\sum_{k=1}^{L}a_{\,k}\arctan\frac{1}{b_{\,k}}
\]
mit $a_{\,k}\in\mathbb Z$ und $b_{\,k}\in\mathbb N$.
Die Argumente sind hierbei Stammbrüche $\frac{1}{n}$. Wir brauchen in obige Reihe nur $m=1$ einsetzen und erhalten sofort
\begin{align*}
&\arctan \left(\frac{1}{n}\right)
=\sum_{k=0}^{\infty}\left(\left(\frac{-4}{\left(n+1\right)^4}\right)^{M}\right)^k
\;\left(\sum_{\ell=0}^{M-1}\left(\frac{-4}{\left(n+1\right)^4}\right)^{\ell}\right.\times
\\
&\Bigg.
\left(
 \frac{1}{\left(4\,k\,M+4\,\ell+1\right)\left(n+1\right)}
+\frac{2}{\left(4\,k\,M+4\,\ell+2\right)\left(n+1\right)^2}
+\frac{2}{\left(4\,k\,M+4\,\ell+3\right)\left(n+1\right)^3}
\right)
\Bigg)
\end{align*}
Die Basis hat sich jetzt vereinfacht zu
\[b=\left(\frac{-4}{\left(n+1\right)^4}\right)^{M}\;.\]
\subsection{1-Term Formel}
Es gibt eine einzige Darstellung mit einem Term, nämlich $\frac{\pi}{4}=\arctan\left(1\right).$
$\\ $
Wir setzen $n=1$ in obige Formel ein und erhalten
\begin{align*}
&\frac{\pi}{4}
=\left(\sum_{k=0}^{\infty}\left(\left(-\frac{1}{4}\right)^M\right)^k\times\right.
\\
&\left(\sum_{\ell=0}^{M-1}\left(-\frac{1}{4}\right)^\ell
\left(
 \frac{1}{2\left({4\,k\,M}+4\,\ell+1\right)}
+\frac{1}{2\left({4\,k\,M}+4\,\ell+2\right)}
+\frac{1}{4\left({4\,k\,M}+4\,\ell+3\right)}
\right)\right)
\end{align*}
In diesem speziellen Fall haben wir die Basis $\left(-\frac{1}{4}\right)^M$ and somit binäre BBP-Formeln für $\frac{\pi}{4}$. Für deren spezielle Eigenschaften verweisen wir auf die Veröffentlichungen
$\\ $
$\\ $
Wir sehen uns die Reihen für einzelne $M$ einmal an:
\[\underline{ M=1}\]
\[\frac{\pi}{4}=\sum_{k=0}^{\infty}
\left(-\frac{1}{4}\right)^k
\left(
 \frac{1}{2\,\left(4\,k+1\right)}
+\frac{1}{2\,\left(4\,k+2\right)}
+\frac{1}{4\,\left(4\,k+3\right)}
\right)
\]
$\\ $
Diese Reihe wurde 1997 von Victor Adamchik and Stan Wagon \cite{AW} veröffentlicht. Sie war bereits Euler bekannt \cite{LE}, wenn auch in etwas andereren Schreibweise (Euler benutzte für die Variablen Planetensymbole).
\[\underline{ M=2}\]
\[\frac{\pi}{4}=\sum_{k=0}^{\infty}\left(\frac{1}{16}\right)^k\times\]
\[
\left(
 \frac{1}{2\,\left(8\,k+1\right)}
+\frac{1}{2\,\left(8\,k+2\right)}
+\frac{1}{4\,\left(8\,k+3\right)}
\right.
\]
\[
\left.
-\frac{1}{8\,\left(8\,k+5\right)}
-\frac{1}{8\,\left(8\,k+6\right)}
-\frac{1}{16\,\left(8\,k+7\right)}
\right)
\]

\[\underline{ M=3}\]
\[\frac{\pi}{4}=\sum_{k=0}^{\infty}\left(-\frac{1}{64}\right)^k\times\]
\[
\left(
 \frac{1}{2\,\left(12\,k+1\right)}
+\frac{1}{2\,\left(12\,k+2\right)}
+\frac{1}{4\,\left(12\,k+3\right)}
\right.
\]
\[
-\frac{1}{8\,\left(12\,k+5\right)}
-\frac{1}{8\,\left(12\,k+6\right)}
-\frac{1}{16\,\left(12\,k+7\right)}
\]
\[
\left.
+\frac{1}{32\,\left(12\,k+ 9\right)}
+\frac{1}{32\,\left(12\,k+10\right)}
+\frac{1}{64\,\left(12\,k+11\right)}
\right)
\]

\[\underline{ M=4}\]
\[\frac{\pi}{4}=\sum_{k=0}^{\infty}\left(\frac{1}{256}\right)^k\times
\]
\[
\left(
 \frac{1}{2\,\left(16\,k+1\right)}
+\frac{1}{2\,\left(16\,k+2\right)}
+\frac{1}{4\,\left(16\,k+3\right)}
\right.
\]
\[
-\frac{1}{8\,\left(16\,k+5\right)}
-\frac{1}{8\,\left(16\,k+6\right)}
-\frac{1}{16\,\left(16\,k+7\right)}
\]
\[
+\frac{1}{32\,\left(16\,k+ 9\right)}
+\frac{1}{32\,\left(16\,k+10\right)}
+\frac{1}{64\,\left(16\,k+11\right)}
\]
\[
\left.
-\frac{1}{128\,\left(16\,k+13\right)}
-\frac{1}{128\,\left(16\,k+14\right)}
-\frac{1}{256\,\left(16\,k+15\right)}
\right)
\]
%
%
\[\vdots \]
\subsection{2-Term Formeln}
Nun kommen wir zu den Formeln mit zwei Summanden. Solche gibt es nach Carl Størmer \cite{CS} genau vier:
\[\frac{\pi}{4}=\arctan\frac{1}{2}+\arctan\frac{1}{3}\]
\[\frac{\pi}{4}=2\arctan\frac{1}{2}-\arctan\frac{1}{7}\]
\[\frac{\pi}{4}=2\arctan\frac{1}{3}+\arctan\frac{1}{7}\]
\[\frac{\pi}{4}=4\arctan\frac{1}{5}-\arctan\frac{1}{239}\]
$ \\ $
Die dritte Formel ist keine unabhängige Formel. Sie lässt sich aus den ersten beiden herleiten. Die Formeln werden verschiedenen Entdeckern zugeschrieben, so Euler, Hermann, Hutton u.a. Die letzte ist die populärste und wurde erstmals von Machin veröffentlicht.
$ \\ $
$ \\ $
Wir beginnen mit der Formel von Euler
\[\frac{\pi}{4}=\arctan\left(\frac{1}{2}\right)+\arctan\left(\frac{1}{3}\right)\]
Wir erhalten die zwei Reihen
\begin{align*}
&\arctan\left(\frac{1}{2}\right)
=\sum_{k=0}^{\infty}\left(\left(-{\frac{4}{81}}\right)^{M}\right)^{k}\times
\\
&\Biggl(\sum_{\ell=0}^{M-1}\left(-{\frac{4}{81}}\right)^{\ell}
\left(
 \frac{1}{3\left(4\,k\,M+4\,\ell+1\right)}
+\frac{2}{9\left(4\,k\,M+4\,\ell+2\right)}
+\frac{2}{27\left(4\,k\,M+4\,\ell+3\right)}
\right)
\Biggr)
\end{align*}
und
\begin{align*}
&\arctan\left(\frac{1}{3}\right)
=\sum_{k=0}^{\infty}\left(\left(-{\frac{1}{64}}\right)^{M}\right)^{k}\times
\\
&\Biggl(\sum_{\ell=0}^{M-1}\left(-{\frac{1}{64}}\right)^{\ell}
\left(
 \frac{1}{ 4\left(4\,k\,M+4\,\ell+1\right)}
+\frac{1}{ 8\left(4\,k\,M+4\,\ell+2\right)}
+\frac{1}{32\left(4\,k\,M+4\,\ell+3\right)}
\right)
\Biggr)
\end{align*}
$ \\ $
Die nächste Formel ist
\[\frac{\pi}{4}=2\arctan\left(\frac{1}{2}\right)-\arctan\left(\frac{1}{7}\right)\]
Die Reihe für $\arctan\left(\frac{1}{2}\right)$ haben wir bereits. Wir benötigen nur noch die Reihe für $\arctan\left(\frac{1}{7}\right)$.
\begin{align*}
&\arctan\left(\frac{1}{7}\right)
=\sum_{k=0}^{\infty}\left(\left(-\frac{1}{1024}\right)^{M}\right)^{k}\times
\\
&\Biggl(\sum_{\ell=0}^{M-1}\left(-\frac{1}{1024}\right)^{\ell}
\left(
 \frac{1}{  8\left(4\,k\,M+4\,\ell+1\right)}
+\frac{1}{ 32\left(4\,k\,M+4\,\ell+2\right)}
+\frac{1}{256\left(4\,k\,M+4\,\ell+3\right)}
\right)
\Biggr)
\end{align*}
Für $M=1$ ergeben sich für  $\arctan\left(\frac{1}{3}\right)$ und $\arctan\left(\frac{1}{7}\right)$ Reihen, die bereits Euler \cite{LE} kannte.
$ \\ $
$ \\ $
Die Reihen für
\[\frac{\pi}{4}=2\arctan\left(\frac{1}{3}\right)+\arctan\left(\frac{1}{7}\right)\]
wurden bereits berechnet. Nun kommen wir zur berühmten Formel von Machin
\[\frac{\pi}{4}=4\arctan\left(\frac{1}{5}\right)-\arctan\left(\frac{1}{239}\right)\]
Wir erhalten die zwei Reihen
\begin{align*}
&\arctan\left(\frac{1}{5}\right)=
\sum_{k=0}^{\infty}\left(\left(-\frac{1}{324}\right)^{M}\right)^{k}\times
\\
&\Biggl(\sum_{\ell=0}^{M-1}\left(-{\frac{1}{324}}\right)^{\ell}
\left(
 \frac{1}{ 6\left(4\,k\,M+4\,\ell+1\right)}
+\frac{1}{18\left(4\,k\,M+4\,\ell+2\right)}
+\frac{1}{108\left(4\,k\,M+4\,\ell+3\right)}
\right)
\Biggr)
\end{align*}
und
\begin{align*}
&\arctan\left(\frac{1}{239}\right)=\sum_{k=0}^{\infty}
\left(\left(-{\frac{1}{829440000}}\right)^{M}\right)^{k}\times
\Biggl(\sum_{\ell=0}^{M-1}\left(-{\frac{1}{829440000}}\right)^{\ell}
\\
&\Biggl.
\left(
 \frac{1}{240\left(4\,k\,M+4\,\ell+1\right)}
+\frac{1}{28800\left(4\,k\,M+4\,\ell+2\right)}
+\frac{1}{6912000\left(4\,k\,M+4\,\ell+3\right)}
\right)
\Biggr)
\end{align*}
\subsection{Formel mit mehr als zwei Termen
}Die nächste Stufe sind jetzt Formeln, die aus mehr als zwei Termen bestehen. Wir geben die effektivsten mit drei und vier Termen an und überlassen die Reihendarstellungen dem Leser:
\[\frac{\pi}{4}
=12\arctan\left(\frac{1}{18}\right)
+8\arctan\left(\frac{1}{57}\right)
-5\arctan\left(\frac{1}{239}\right)
\]
\[\frac{\pi}{4}
=44\arctan\left(\frac{1}{57}\right)
+7\arctan\left(\frac{1}{239}\right)
-12\arctan\left(\frac{1}{682}\right)
+24\arctan\left(\frac{1}{12943}\right)
\] 
Für weitere $\frac{\pi}{4}$-Darstellungen vom Machin'schen Typ verweisen wir auf die sehr ausführliche Darstellung von Jörg Arndt \cite{JA}. Je mehr Terme verwendet werden, umso kleiner wird das Argument und umso kürzer ist die Rechenzeit. Ist zum Beispiel
$\arctan\left(\frac{1}{n}\right)$ zu berechnen mit einem schon sehr kleinen Argument $\frac{1}{n}$, so führt die Verwendung der hier vorgestellten Methode zu einem noch kleinerem Argument der Größenordnung
$\left(\frac{-4}{(n+1)^4}\right)^{M}.$

\end{document}